\input epsf
\documentclass[a4paper,11pt]{amsart}
\usepackage{latexsym}
\newtheorem{thm}{Theorem}[section]
\newtheorem{prop}[thm]{Proposition}

\title{Le cocycle du verger}
\author{Roland Bacher}

\begin{document}

\maketitle         
 
{\begin{abstract}
We define and prove uniqueness of a 
natural homomorphism (called the Orchard morphism)
from some groups associated naturally to a finite set $E$ to the
group ${\mathcal E}(E)$ of two-partitions of $E$ representing 
equivalence relations having at most two classes on $E$.

As an application, given a finite generic configuration ${\mathcal C}\subset
{\mathbf R}^d$, we exhibit a natural partition of ${\mathcal C}$ in two sets.

Le but de cette note est de montrer l'existence et l'unicit\'e d'un 
homomorphisme naturel non-trivial
entre certains groupes associ\'es \`a un ensemble fini.

Cet homomorphisme fournit une partition naturelle en deux 
sous-ensembles sur l'ensemble ${\mathcal C}\subset
{\mathbf R}^d$ des points d'une configuration finie g\'en\'erique.
\end{abstract} }\footnotetext{Keywords: Group, Configuration of points, 
Two-partition, Semi-orientation, Orchard morphism. AMS-Classification: 
05C25, 52C35}

Soit $E$ un ensemble fini.
Notons $E^{(l)}$ l'ensemble des suites de longueur $l$ sans r\'ep\'etitions
et ${\mathcal F}_+(E^{(l)})$ l'ensemble des fonctions de $E^{(l)}$
dans $\{\pm 1\}$ qui sont sym\'etriques, c'est-\`a-dire ind\'ependantes de
l'ordre des arguments. Notons $\Delta_E$ le complexe simplicial fini
associ\'e au simplexe de sommets $E$. En utilisant l'isomorphisme
entre le groupe additif 
${\mathbf Z}/2{\mathbf Z}$ et le groupe multiplicatif $\{\pm 1\}$,
le groupe multiplicatif ${\mathcal F}_+(E^{(l)})$
s'identifie au groupe $C^{l-1}$ des $(l-1)-$cocha\^{\i}nes
simpliciales de $\Delta_E$ \`a valeurs dans ${\mathbf Z}/2{\mathbf Z}$.
L'homomorphisme cobord $\partial:C^{l-1}\longrightarrow C^l$
de ${\mathcal F}_+(E^{(l)})$ dans 
${\mathcal F}_+(E^{(l+1)})$ est donn\'ee par
$$\partial\varphi(x_0,\dots,x_l)=\prod_{j=0}^l\varphi(x_0,\dots,
x_{j-1},\hat x_j,x_{j+1},\dots,x_l)$$
et on a \'evidemment les inclusion
$$B^k=\hbox{Im}\{\partial:C^{k-1}\rightarrow C^{k}\}\subset
Z^k=\hbox{Ker}\{\partial:C^{k}\rightarrow C^{k+1}\}\subset
C^k={\mathcal F}_+(E^{(k+1)})$$
du groupe $B^k$ des $k-$cobords dans le groupe $Z^k$ des $k-$cocycles
et du groupe $Z^k$
dans le groupe $C^k={\mathcal F}(E^{(k+1)})$ des $k-$cocha\^{\i}nes.

Comme $\Delta_E$ est homotope \`a un point, les groupes de cohomologie
$H^i=Z^i/B^i$
sont triviaux pour $i\geq 1$ et $H^{0}=\{\pm 1\}$ (cf. \cite{Br} ou \cite{M}).
Le sous-groupe $Z^1\subset {\mathcal F}_+(E^{(2)})$
des $1-$cocycles coincide donc avec le groupe $B^1$
des $1-$cobords et s'identifie au groupe multiplicatif
${\mathcal E}(E)=C^0/Z^0$ dont les \'el\'ements sont des 
paires de fonctions $\{\pm \alpha\}$
de $E$ dans $\{\pm 1\}$. Une telle paire de fonctions $\pm \alpha$
d\'efinit une partition $\alpha^{-1}(1)\cup\alpha^{-1}(-1)$ de $E$ en 
au plus deux sous-ensembles non-vides et correspond au $1-$cobord
$\sigma(x,y)=\partial\alpha(x,y)=\alpha(x)\alpha(y)$ pour
$x\not= y\in E$. 

R\'eciproquement, une cocha\^{\i}ne 
$\sigma\in C^1={\mathcal F}_+(E^{(2)})$ est un $1-$cocycle et donc un 
$1-$cobord
si $\partial\sigma(a,b,c)=\sigma(a,b)\ \sigma(a,c)\ \sigma(b,c)=1$
pour tout $(a,b,c)\in E^{(3)}$. 
On construit alors $\alpha\in\partial^{-1}(\sigma)\subset C^0={\mathcal F}_+
(E^{(1)})$ en 
choisissant un point base $x_0\in E$ et en posant
$\alpha(x_0)=1$ et $\alpha(x)=\sigma(x,x_0)$ pour $x\not=x_0$.
(L'existence de la cocha\^{\i}ne $\alpha\in\partial^{-1}\sigma$ peut
\'egalement se d\'eduire des propri\'et\'es du graphe fini $\Gamma$ 
de sommets $E$ et d'ar\^etes $\{x,y\}$ pour $\sigma(x,y)=1$. 
L'identit\'e $\sigma(a,b)\ \sigma(a,c)\ \sigma(b,c)=1$ 
pour tout $(a,b,c)\in E^{(3)}$ montre
que $\Gamma$ poss\`ede au plus deux composantes
connexes $\Gamma_+$ et $\Gamma_-$ qui sont des graphes complets.
Quitte \`a permuter $\Gamma_+$ et $\Gamma_-$ on peut supposer $x_0\in
\Gamma_+$ et on d\'efinit $\alpha(x)=1$ si $x\in\Gamma_+$ et $\alpha(x)=-1$
sinon.)
J'appellerai ${\mathcal E}(E)=\{\pm 1\}^E/\pm 1$ 
le groupe des {\it $2-$partitions}
de $E$. Ses \'el\'ements sont en bijection avec les relations
d'\'equivalence sur $E$ ayant au plus $2$ classes. 

Introduisons maintenant le groupe multiplicatif
$${\mathcal F}_\pm(E^{(l)})={\mathcal F}_+(E^{(l)})\cup
{\mathcal F}_-(E^{(l)})$$
des fonctions de $E^{(l)}$ dans $\{\pm 1\}$ qui sont soit
sym\'etriques soit antisym\'etriques. Une telle fonction $\varphi
\in {\mathcal F}_\epsilon(E^{(l)})$ (avec $\epsilon=+$ ou $\epsilon=-$)
v\'erifie l'\'egalit\'e
$$\varphi(x_1,\dots,x_{i+1},x_i,\dots,x_l)=\epsilon\ \varphi(x_1,\dots
,x_i,x_{i+1},\dots,x_l)$$
pour tout $1\leq i<l$.

Munissons l'ensemble fini $E$ d'un ordre total. Pour un sous-ensemble
$F\subset E$, d\'esignons par 
${F\choose k}$ l'ensemble des suites strictement croissantes de longueur
$k$ dans $F$. Notons $\sharp(E)$ le nombre d'\'el\'ements de $E$ et posons
pour $\varphi\in{\mathcal F}_\pm (E^{(l)})$ et $y\not= z\in E$
$$\sigma_\varphi(y,z)=\epsilon^{\sharp(E)-3\choose l-2}
\prod_{(x_1,\dots,x_{l-1})\in{E\setminus\{y,z\}\choose l-1}}
\varphi(x_1,\dots,x_{l-1},y)\varphi(x_1,\dots,x_{l-1},z)$$
o\`u $\epsilon=1$ si $\varphi\in{\mathcal F}_+(E^{(l)})$ et
$\epsilon=-1$ si $\varphi\in{\mathcal F}_-(E^{(l)})$.

\begin{prop} \label{cocycleverger}
(i) La fonction $\sigma_\varphi:E^{(2)}\longrightarrow \{\pm 1\}$ 
est sym\'etrique et ne d\'epend pas du choix de l'ordre sur $E$.

\ \ (ii) L'application $\varphi\longmapsto \sigma_\varphi$ d\'efinit un 
homomorphisme $\hbox{Sym}(E)-$\'equivariant du groupe ${\mathcal F}_\pm
(E^{(l)})$ dans le sous-groupe $Z^1\subset{\mathcal F}_+(E^{(2)})$
des $1-$cocycles.
\end{prop}

\noindent{\bf Esquisse de la preuve.} La sym\'etrie de $\sigma_\varphi$ est \'evidente. 
Le choix d'un autre ordre de $E$ permute de mani\`ere similaire
les arguments dans les deux facteurs 
$$\varphi(x_1,\dots,x_{l-1},y)\varphi(x_1,\dots,x_{l-1},z)$$ 
et les changements \'eventuels de signe s'annulent deux-\`a-deux.

Les d\'etails pour prouver que
l'application $\varphi\longmapsto \sigma_\varphi$ 
est un homomorphisme $\hbox{Sym}(E)-$\'equivariant
sont faciles et laiss\'es au lecteur. 

Le cocycle $\sigma_\varphi$ est ferm\'e si on a 
$\sigma_\varphi(a,b)\sigma_\varphi(b,c)\sigma_\varphi(a,c)=1$
pour tout $(a,b,c)\in E^{(3)}$.
Si $\varphi$ est sym\'etrique cette identit\'e est facile et laiss\'ee 
au lecteur. Pour $\varphi$ antisym\'etrique, elle r\'esulte
de l'\'egalit\'e
$$\begin{array}{l}
\displaystyle \varphi(x_1,\dots,x_{l-2},c,a)\ 
\varphi(x_1,\dots,x_{l-2},c,b) \\  
\displaystyle \varphi(x_1,\dots,x_{l-2},a,b)\ 
\varphi(x_1,\dots,x_{l-2},a,c) \\  
\displaystyle \varphi(x_1,\dots,x_{l-2},b,a)\ 
\varphi(x_1,\dots,x_{l-2},b,c)=-1\end{array}$$
pour $(x_1,\dots,x_{l-2})\in{E\setminus \{a,b,c\}\choose l-2}$.\hfill
$\Box$

\smallskip

On appelle le cocycle $\sigma_\varphi$ donn\'e par la Proposition
\ref{cocycleverger} le {\it cocycle du verger}.
Le {\it morphisme du verger} est d\'efini comme l'homomorphisme
de ${\mathcal F}_\pm (E^{(l)})$ dans ${\mathcal E}(E)$
qui associe \`a la fonction sym\'etrique ou antisym\'etrique
$\varphi\in {\mathcal F}_\pm (E^{(l)})$ la $2-$partition
$\partial^{-1}(\sigma_\varphi)$ 
d\'efinie par le $1-$cocycle $\sigma_\varphi$.

Pour $1\leq l<\sharp(E)$, on peut montrer que c'est 
l'unique homomorphisme non-trivial du groupe
${\mathcal F}_\pm (E^{(l)})$ dans le groupe ${\mathcal E}(E)$ qui
soit \'equivariant par rapport \`a l'action \'evidente par
automorphismes de $\hbox{Sym}(E)$ sur ${\mathcal F}_\pm (E^{(l)})$
et ${\mathcal E}(E)$.

Pour $l=\sharp(E)$ le morphisme du verger est toujours trivial
et il n'existe pas d'autre homomorphisme naturel 
($\hbox{Sym}(E)-$\'equivariant)
de ${\mathcal F}_\pm(E^{(\sharp(E))})$ dans ${\mathcal E}(E)$ 
si $\sharp(E)>2$. Pour $\sharp(E)=2$, 
un homomorphisme naturel \lq\lq exotique'' (distinct du morphisme du verger
qui est trivial dans ce cas) de 
${\mathcal F}_\pm(E^{(\sharp(E))})$ 
dans ${\mathcal E}(E)$ existe cependant ; c'est l'application qui envoie 
une fonction antisym\'etrique (qui est unique, au signe pr\`es)
de ${\mathcal F}_-(E^{(2)})$ 
sur l'unique \'el\'ement non-trivial de ${\mathcal E}(E)$.  
Cette exception est rendue possible par le fait que le groupe
$\hbox{Sym}(E)$ agit trivialement sur ${\mathcal E}(E)$ si $E$ 
poss\`ede au plus deux \'el\'ements.

L'existence de cet homomorphisme n'est pas surprenant sur
les fonctions sym\'etriques. Il peut en effet s'obtenir
en composant l'un des homomorphismes
$$\begin{array}{lcl}\displaystyle 
\varphi&\longmapsto& \beta_\varphi(y)=\prod_{(x_1,\dots,x_{l-1})\in{E\setminus\{y\}\choose l-1}}
\varphi(x_1,\dots,x_{l-1},y)\\
\displaystyle \varphi&\longmapsto &
\tilde \beta_\varphi(y)=\prod_{(x_1,\dots,x_{l})\in{E\setminus\{y\}\choose l}}
\varphi(x_1,\dots,x_{l})\end{array}$$
de ${\mathcal F}_+(E^{(l)})\longrightarrow {\mathcal F}_+(E^{(1)})$ 
avec l'homomorphisme quotient ${\mathcal F}_+(E^{(1)})\longrightarrow
{\mathcal F}_+(E^{(1)})/\{\pm 1\}={\mathcal E}(E)$.
C'est son extension
aux fonctions antisym\'etriques ${\mathcal F}_-(E^{(l)})$ qui n'est pas
\'evidente \`a priori et qui constitue le r\'esultat
de cette Note.

\smallskip

Citons l'application suivante, \`a l'origine du nom 
de cet homomorphisme.

Une {\it configuration} est un sous-ensemble fini ${\mathcal C}\subset 
{\mathbf R}^d$
de points dans un espace vectoriel (ou affine) r\'eel muni
d'une orientation. Une telle configuration est 
{\it g\'en\'erique} si tout sous-ensemble
de $k+1\leq d+1$ points engendre un sous-espace affine de dimension $k$.
Consid\'erons la fonction
$\varphi:{\mathcal C}^{(d+1)}\longrightarrow\{\pm 1\}$
d\'efinie par $\varphi(x_0,\dots,x_d)=1$ si la base
$x_1-x_0,x_2-x_1,,\dots,x_d-x_{d-1}$ induit l'orientation positive de
${\mathbf R}^d$ et $\varphi(x_0,\dots,x_d)=-1$ sinon. Il est \'evident
que $\varphi$ est une fonction antisym\'etrique. On
obtient alors une $2-$partition sur ${\mathcal C}$ en 
consid\'erant l'image $\partial^{-1}(\sigma_\varphi)$
de $\varphi$ par le morphisme du verger.

Ainsi, dans un chapitre apocryphe d'Alice au pays des Merveilles,
la Reine de Coeur exige un verger plant\'e
de fa\c con rationelle de pruniers et cerisiers. Elle a impos\'e
les emplacements des arbres futurs selon une configuration
g\'en\'erique, semant ainsi un grand
d\'esarroi parmi ses jardiniers. Alice les a sauv\'es
provisoirement de la d\'ecapitation en utilisant le
morphisme du verger qui g\'en\'eralise aux dimensions
sup\'erieures les haies plant\'ees de deux essences en
alternance.

La figure montre le verger de la Reine, plant\'e par
Alice avec $3+6=9$ arbres (Alice a quand-m\^eme eu un moment 
d'h\'esitation : ne sachant pas si la Reine pr\'ef\'erait les
prunes ou les cerises, elle a jou\'e \`a pile ou face). 
Les droites trac\'ees ont \'et\'e utilis\'ees pour convaincre
la Reine (qui n'y a rien compris mais qui, craignant de perdre la 
face devant ses sujets, n'a pas os\'e l'avouer)
en utilisant les explications donn\'ees ci-dessous.
\vskip0.5cm

\centerline{\epsfysize8cm\epsfbox{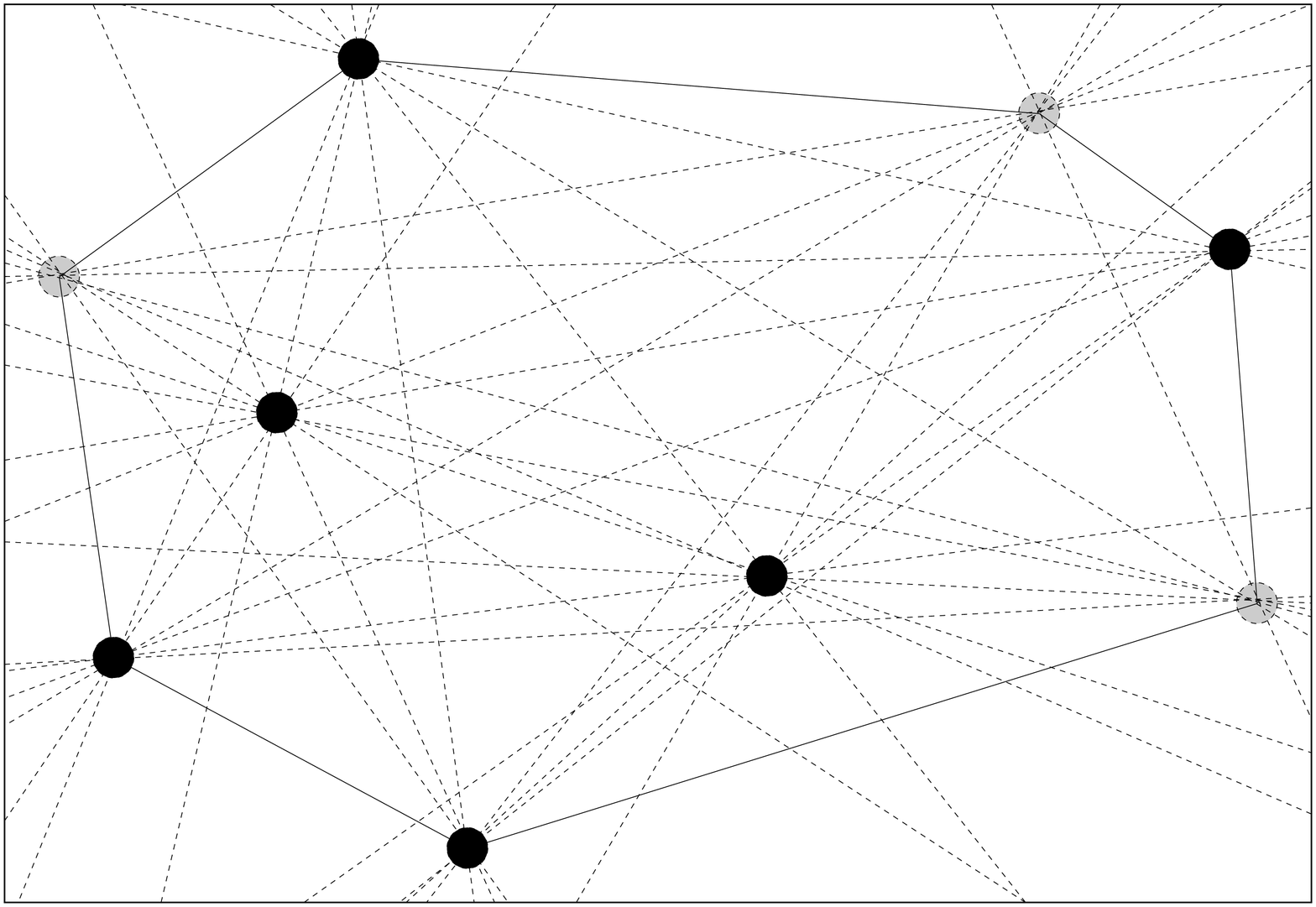}}
\centerline{Figure : Le verger de la Reine de Coeur} 
\medskip

Pour $\varphi\in{\mathcal F}_-({\mathcal C}^{d+1})$
provenant d'une configuration g\'en\'erique ${\mathcal C}\subset 
{\mathbf R}^d$,
la relation d'\'equivalence associ\'ee \`a $\sigma_\varphi$ se calcule
de la mani\`ere suivante: Pour $x\not= y\in{\mathcal C}$, notons $s(x,y)$
le nombre d'hyperplans qui s\'eparent $x$ de $y$ et qui contiennent 
$d$ points distincts de ${\mathcal C}\setminus\{x,y\}$.
Les points $x$ et $y$ appartiennent
\`a la m\^eme classe de $\partial^{-1}(\sigma_\varphi)$
si et seulement si $s(x,y)\equiv {\sharp({\mathcal C})-3
\choose d-1}\pmod 2$.

Citons pour terminer la propri\'et\'e suivante de 
l'homomorphisme du verger. On dit que deux fonctions $\varphi,\psi\in
{\mathcal F}_\pm(E^{(l)})$ sont {\it reli\'ees par un 
flip} si 
$$\varphi(x_1,\dots,x_l)=\psi(x_1,\dots,x_l)$$
sauf pour $\{x_1,\dots,x_l\}=X$ o\`u $X\subset E$
est un sous-ensemble fix\'e de $l$ \'el\'ements dans $E$. 
Les $2-$partitions $\partial^{-1}(\sigma_\varphi)$
et $\partial^{-1}(\sigma_\psi)$ associ\'ees \`a $\varphi$ et $\psi$
sont alors identiques en dehors de $X$ et diff\`erent sur $X$ 
(pour un choix convenable des repr\'esentants).

Sur une configuration planaire, 
un flip consiste \`a faire d\'eg\'en\'erer un triangle tr\`es
obtus form\'e de trois points convenables de ${\mathcal C}$ en trois points
align\'es pour arriver \`a un triangle tr\`es obtus ressemblant
au triangle mirroir relativement \`a son ar\^ete la plus longue. En coloriant
les deux classes associ\'ees aux $2-$partitions 
avec deux couleurs diff\'erentes,
on voit ainsi que lors d'un flip, tous les points restent de la 
m\^eme couleur sauf les points du flip qui, sous l'effet
de la forte \'emotion provoqu\'ee par l'alignement,
changent de couleur. 

Remarquons que deux configurations g\'en\'eriques de $n$ points
dans ${\mathbf R}^2$ (ou plus g\'en\'eralement dans ${\mathbf R}^d$) 
peuvent toujours \^etre reli\'ees, \`a isotopie pr\`es, par un nombre fini
de flips.   

Des calculs indiqu\'es dans \cite{BG} sugg\`erent qu'il n'y a g\'en\'eralement
pas de restrictions sur les nombres possibles de pruniers 
dans les configurations planaires g\'en\'eriques
\`a $n$ points : Pour $n=5,7,8,9$ ce nombre peut prendre toutes les 
valeurs entre $0$ et $n$. Les cas monochromatiques ou tr\`es
d\'es\'equilibr\'es sont cependant relativement rares.

J'aimerais remercier toutes les personnes avec qui j'ai eu des 
discussions sur les vergers, en particulier M. Brion, P. Cameron, N. A'Campo,
E. Ferrand, E. Ghys, P. de la Harpe et A. Marin.

\vskip1cm

Roland Bacher

INSTITUT FOURIER

Laboratoire de Math\'ematiques

UMR 5582 (UJF-CNRS)

BP 74

38402 St MARTIN D'H\`ERES Cedex (France)
 
e-mail: Roland.Bacher@ujf-grenoble.fr  


\begin{thebibliography}{}

\bibitem{B} Bacher R, {\it An Orchard Theorem}, Preprint math/CO0206266.

\bibitem{BG} Bacher R, Garber D, {\it Chromatic properties of generic planar
configurations of points}, Preprint math/GT0210051.

\bibitem{Br} Bredon G.E, {\it Topology and Geometry}, Springer (1993).

\bibitem{M} Massey W.S, {\it A Basic Course in Algebraic Topology},
Springer (1991).

\end{thebibliography}
\end{document}